\documentclass[a4paper,12pt]{article}

\usepackage{amsthm}
\usepackage{amsmath,amssymb,latexsym,amsfonts,mathrsfs}

\renewcommand{\d}{{\rm d}} 
\newcommand{\sgn}{\mathop{\rm sgn}}
\newcommand{\eps}{\ensuremath{\varepsilon}}
\newcommand{\e}{{\rm e}} 

\newcommand{\tend}[2]{\mathrel{\mathop{\longrightarrow}\limits^{#1}_{#2}}}

\renewcommand{\bar}{\overline}

\newcommand{\down}{\ensuremath{\downarrow}}

\newcommand{\rbra}[1]{\!\left( #1 \right)} 
\newcommand{\cbra}[1]{\!\left\{ #1 \right\}} 
\newcommand{\sbra}[1]{\!\left[ #1 \right]} 

\newcommand{\bD}{\ensuremath{\mathbb{D}}}
\newcommand{\bN}{\ensuremath{\mathbb{N}}}
\newcommand{\bR}{\ensuremath{\mathbb{R}}}

\newcommand{\cB}{\ensuremath{\mathcal{B}}}
\newcommand{\cE}{\ensuremath{\mathcal{E}}}
\newcommand{\cF}{\ensuremath{\mathcal{F}}}
\newcommand{\cG}{\ensuremath{\mathcal{G}}}

\newcommand{\sP}{\ensuremath{\mathscr{P}}}
\newcommand{\sW}{\ensuremath{\mathscr{W}}}

\newcommand{\ve}{\ensuremath{\boldsymbol{e}}}

\theoremstyle{plain}
\newtheorem{Thm}{Theorem}[section]

\newtheorem{Prop}[Thm]{Proposition}

\theoremstyle{definition}

\newcommand{\Proof}[2][Proof]{\begin{proof}[{#1}] #2 \end{proof}}

\numberwithin{equation}{section}

\makeatletter
\renewcommand\section{\@startsection {section}{1}{\z@}%
                                   {-3.5ex \@plus -1ex \@minus -.2ex}%
                                   {2.3ex \@plus.2ex}%
                                   {\normalfont\large\bf}}
\makeatother

\makeatletter
\renewcommand\subsection{\@startsection {subsection}{1}{\z@}%
                                   {-3.5ex \@plus -1ex \@minus -.2ex}%
                                   {2.3ex \@plus.2ex}%
                                   {\normalfont\normalsize\bf}}
\makeatother

\begin{document}
\begin{center}
{\Large \bf On universality in penalisation problems with multiplicative weights} 
\end{center}

\begin{center}
Kouji Yano\footnote{
Graduate School of Science, Kyoto University.}\footnote{
The research of Kouji Yano was supported by
JSPS KAKENHI grant no.'s JP19H01791, JP19K21834, JP21H01002 and JP18K03441 
and by JSPS Open Partnership Joint Research Projects grant no. JPJSBP120209921.}
\end{center}

\begin{abstract}
We give a general framework for the universality classes of $ \sigma $-finite measures 
in penalisation problems with multiplicative weights. 
We discuss penalisation problems for Brownian motions, L\'evy processes and Langevin processes 
in our framework. 
\end{abstract}


\section{Introduction}

For a measure $ \mu $ and a non-negative measurable function $ f $, 
we write $ \mu[f] $ for the integral $ \int f \d \mu $. 

For a probability space $ (\Omega,\cF,P) $ 
equipped with a filtration $ (\cF_s)_{s \ge 0} $, 
and for a non-negative process $ \Gamma = (\Gamma_t)_{t \ge 0} $ called a \emph{weight}, 
we mean by a {\em penalisation} a problem 
of finding a limit probability $ P^{\Gamma} $ on $ (\Omega,\cF) $ 
called the \emph{penalised probability} such that 
\begin{align}
\frac{P[F_s \Gamma_t]}{P[\Gamma_t]} \tend{}{t \to \infty } P^{\Gamma}[F_s] 
\label{}
\end{align}
is satisfied 
for all $ s \ge 0 $ and all bounded $ \cF_s $-measurable functional $ F_s $. 
Under the penalised probability $ P^{\Gamma} $, 
the process $ (\Gamma_t)_{t \ge 0} $ 
is prevented from taking small values; 
this is why Roynette--Vallois--Yor \cite{MR2261065} (see also \cite{RoyYor}) 
called this problem the penalisation. 
Conditioning a process to stay in a domain $ D $ 
may be regarded as a special case of the penalisation, 
as we take the weight 
$ \Gamma_t = 1_{\{ \tau_D > t \}} $ 
where $ \tau_D $ denotes the exit time of $ D $.

Although the penalised probability $ P^{\Gamma} $ depends upon the weight $ \Gamma $, 
we can often find a $ \sigma $-finite measure $ \sP $ on $ (\Omega,\cF) $ 
independent of a particular weight such that 
\begin{align}
P^{\Gamma}(A) = \frac{\sP[\Gamma_{\infty };A]}{\sP[\Gamma_{\infty }]} 
, \quad A \in \cF 
\label{}
\end{align}
holds with a suitable limit $ \Gamma_{\infty } $ of $ \Gamma_t $ 
in a certain class of weights $ \Gamma $. 
In this case we say that $ \Gamma $ belongs to 
the \emph{universality class} of $ \sP $. 
The aim of this paper is to gain a clear insight into the universality classes 
in penalisation problems.  
For this purpose, we confine ourselves to multiplicative weights. 

Let $ \{ B=(B_t)_{t \ge 0},W_x \} $ denote the canonical representation 
of the one-dimensional Brownian motion with $ W_x(B_0=x)=1 $ 
and let $ \cF^B_t = \sigma(B_s : s \le t) $ denote the natural filtration 
of the coordinate process $ B $. 
Let $ \tau_D = \inf \{ t \ge 0 : B_t = 0 \} $ denote the exit time of $ B $ 
from the non-zero real $ D = \bR \setminus \{ 0 \} $. 
Let $ x \in D $ be fixed. 
It is then well-known that 
\begin{align}
W_x[F_s|\tau_D>t] \tend{}{t \to \infty } W^{\rm \pm 3B}_x[F_s] 
= \frac{1}{|x|} W_x \sbra{ \Big. F_s |B_s| 1_{\{ \tau_D>s \}} } 
\label{}
\end{align}
for all bounded $ \cF^B_s $-measurable functional $ F_s $, 
where $ W^{\rm \pm 3B}_x $ denotes the law of $ \pm $ times 3-dimensional Bessel process 
starting from $ x $. 
This conditioning to avoid zero may be regarded as a special case of the penalisation 
with the weight being given by $ \Gamma_t = 1_{\{ \tau_D>t \}} $. 
Note that $ W^{\rm \pm 3B}_x $ is locally absolutely continuous with respect to $ W_x $, i.e. 
$ W^{\rm \pm 3B}_x|_{\cF^B_s} $ is absolutely continuous with respect to $ W_x|_{\cF^B_s} $ 
for all $ s \ge 0 $. 
But $ W^{\rm \pm 3B}_x $ and $ W_x $ are mutually singular on $ \cF^B_{\infty } := \sigma(B) $, 
because $ W^{\rm \pm 3B}_x(\tau_D = \infty ) = W_x(\tau_D < \infty ) = 1 $. 
While the original process $ \{ B,W_x \} $ is recurrent, 
the \emph{penalised process} $ \{ B,W^{\rm \pm 3B}_x \} $ is transient.

Roynette--Vallois--Yor 
(\cite{MR2229621} and \cite{MR2253307}) 
have studied the penalisation problems for the one-dimensional Brownian motion. 
They determined the penailsed probabilities 
for $ \Gamma_t = f(\bar{X}_t) $, a function of a supremum, 
$ \Gamma_t = f(L_t) $, a function of a local time at 0, 
and $ \Gamma_t = \exp (- \int_0^t v(B_s) \d s) $, a Kac killing weight. 
For the special case $ \Gamma_t = \e^{- L_t} $, we have 
\begin{align}
\frac{W_0[F_s \e^{-L_t}]}{W_0[\e^{-L_t}]} 
\tend{}{t \to \infty } 
W^{\Gamma}_0[F_s] = \frac{1}{1+|x|} W_0 \sbra{ \Big. F_s(1+|B_s|) \e^{-L_s} } 
\label{}
\end{align}
for all $ s \ge 0 $ and all bounded $ \cF^B_s $-measurable functional $ F_s $. 
Although $ W^{\Gamma}_0 $ is locally absolutely continuous with respect to $ W_0 $, 
the two measures $ W^{\Gamma}_0 $ and $ W_0 $ are mutually singular on $ \cF^B_{\infty } $, 
because $ W^{\Gamma}_0(L_{\infty } < \infty ) = W_0(L_{\infty } = \infty ) = 1 $. 
While the original process $ \{ B,W_0 \} $ is recurrent, 
the penalised process $ \{ B,W^{\Gamma}_0 \} $ is transient.

Najnudel--Roynette--Yor (\cite{NRY}) have introduced 
the $ \sigma $-finite measure $ \sW_0 $ defined by 
\begin{align}
\sW_0 = \int_0^{\infty } \frac{\d u}{\sqrt{2 \pi u}} \, \Pi^{(u)} \bullet W^{\rm s3B}_0 , 
\label{}
\end{align}
where $ \Pi^{(u)} $ stands for the law of the Brownian bridge from 0 to 0 of length $ u $, 
$ W^{\rm s3B}_0  $ for the law of the symmetrised Bessel process, 
and $ \bullet $ for the law of the concatenated path of two independent paths. 
They proved that the penalised probability $ W^{\Gamma}_0 $ 
for any weight $ \Gamma $ in the previous paragraph 
is absolutely continuous on $ \cF^B_{\infty } $ with respect to $ \sW_0 $: 
\begin{align}
W^{\Gamma}_0[F] = \frac{\sW_0[F \Gamma_{\infty }]}{\sW_0[\Gamma_{\infty }]} 
\label{}
\end{align}
for all bounded $ \cF^B_{\infty } $-measurable functional $ F $. 
Moreover, if we define $ \sW_x(\cdot) = \sW_0(x+B \in \cdot) $, we have 
\begin{align}
W^{\rm \pm 3B}_x[F] = \frac{\sW_x[F ; \tau_D = \infty ]}{\sW_x(\tau_D = \infty )} 
\label{}
\end{align}
for all $ x>0 $ 
and all bounded $ \cF^B_{\infty } $-measurable functional $ F $. 
In other words, all the weights belong to the universality calss of $ \sW_x $.

K.Yano--Y.Yano--Yor \cite{MR2552915,MR2744885}, 
Y.Yano \cite{MR3127911} 
and recently Takeda--K.Yano \cite{TY}
studied the penalisation problems 
for one-dimensional stable L\'evy processes 
and found out that there are two different universality classes. 
In this paper, we would like to give a general framework 
to characterise universality classes, 
where we will give some new results.

Groeneboom--Jongbloed--Wellner \cite{MR1733148} 
studied the conditioning to stay positive 
for the Langevin process. 
Profeta \cite{MR3386368} studied penalisation problems with several kinds of weights. 
In this paper, we shall discuss universality classes 
for those penalisation problems.

This paper is organized as follows. 
In Section \ref{sec: penal meas} 
we develop a general study on penalised probabilities 
with multiplicative weights. 
In Section \ref{sec: subs M} 
we define the unweighted measures 
and discuss the subsequent Markov property of them. 
In Section \ref{sec: univ} 
we state and prove our main theorems on universality classes. 
In Section \ref{sec: penal}
we give a general discussion on penalisation problems 
with multiplicative weights. 
In Sections \ref{sec: BM}, \ref{sec: Levy} and \ref{sec: Lan}, 
we look at some known results of penalisation problems 
for Brownian motions, L\'evy processes and Langevin processes in our framework. 
In Section \ref{sec: ext} as an appendix, 
we discuss extension of the transformed probability measures 
given by local absolute continuity.

\section{Penalised probability} \label{sec: penal meas}

For a measure $ \mu $ and a non-negative measurable function $ f $, 
we write $ f \cdot \mu $ for the transformed measure defined by 
$ (f \cdot \mu)(A) = \int_A f \d \mu $ for all measurable set $ A $. 
Let $ (\cF_s)_{s \ge 0} $ be a filtration. 
For two measures $ \mu $ and $ \nu $, 
we say that $ \mu $ is \emph{locally absolutely continuous} with respect to $ \nu $ 
if $ \mu|_{\cF_s} $ is absolutely continuous with respect to $ \nu $. 
We say the two measures are \emph{locally equivalent} 
if they are locally absolutely continuous with respect to each other. 
For a parameterised family $ (\mu_\lambda)_{\lambda} $ of finite measures 
and a finite measure $ \mu $, we say that 
\begin{align}
\text{$ \lim \nolimits_{\lambda} \mu_{\lambda} = \mu $ \emph{along} $ (\cF_s)_{s \ge 0} $} 
\label{}
\end{align}
if 
\begin{align}
\lim \nolimits_{\lambda} \mu_{\lambda}[F_s] = \mu[F_s] 
\label{}
\end{align}
holds for all $ s \ge 0 $ and all bounded measurable functional $ F_s $.

Let $ S $ be a locally compact separable metric space 
and let $ \bD $ denote the space of c\`adl\`ag paths 
from $ [0,\infty ) $ to $ S $. 
Let $ X = (X_t)_{t \ge 0} $ denote the coordinate process: 
$ X_t(\omega) = \omega(t) $ for $ t \ge 0 $ and $ \omega \in \bD $. 
Let $ \cF^X_t = \sigma(X_s:s \le t) $ denote the natural filtration of $ X $
and set $ \cF_t = \bigcap_{\eps>0} \cF^X_{t+\eps} $ so that 
$ (\cF_t)_{t \ge 0} $ is a right-continuous filtration. 
We write $ \cF_{\infty } = \sigma(\bigcup_{t \ge 0} \cF_t) = \sigma(X) $. 
For $ t \ge 0 $, let $ \theta_t $ denote the shift operator of $ \bD $: 
$ \theta_t \omega (s) = \omega(t+s) $ for $ s \ge 0 $.

Let $ \{ X,\cF_{\infty },(P_x)_{x \in S} \} $ denote the canonical representation 
of a strong Markov process taking values in $ S $ 
with respect to the augmented filtration $ (\cG_t)_{t \ge 0} $ of $ (\cF_t)_{t \ge 0} $. 
A process $ \Gamma = (\Gamma_t)_{t \ge 0} $ is called a {\em weight} 
if it is a non-negative c\`adl\`ag process. 
A weight $ \Gamma $ is called {\em multiplicative} if 
$ \Gamma $ is adapted to $ (\cF_t)_{t \ge 0} $ and 
\begin{align}
\text{$ \Gamma_t = \Gamma_s \cdot (\Gamma_{t-s} \circ \theta_s) $, 
$ P_x $-a.s. for all $ 0 \le s \le t < \infty $ and all $ x \in S $.} 
\label{}
\end{align}
Let $ \Gamma $ be a multiplicative weight. 
Since $ \Gamma_0 = \Gamma_0 \cdot (\Gamma_0 \circ \theta_0) = \Gamma_0^2 $, 
we note that 
\begin{align}
\text{for any $ x \in S $ we have either $ P_x(\Gamma_0 = 1)=1 $ or $ P_x(\Gamma_0 = 0)=1 $}. 
\label{}
\end{align}
We set 
\begin{align}
S^{\Gamma} = \cbra{ x \in S : P_x(\Gamma_0=1)=1 } . 
\label{}
\end{align}
It is easy to see that 
\begin{align}
\tau^{\Gamma} 
:= \inf \{ t \ge 0 : X_t \notin S^{\Gamma} \} 
= \inf \{ t \ge 0 : \Gamma_t=0 \} 
\ \text{$ P_x $-a.s. for all $ x \in S $}, 
\label{}
\end{align}
since [$ \Gamma_{t_0} = 0 $ implies $ \Gamma_t=0 $ for all $ t \ge t_0 $] 
because of the multiplicativity. 

We introduce the following assumptions: 
\begin{enumerate}
\item[{\bf (A1)}] 
There is a Borel function $ \varphi^{\Gamma} $ on $ S $ 
such that 
$ \varphi^{\Gamma} > 0 $ on $ S^{\Gamma} $ and 
\begin{align}
P_x[ \Gamma_t \varphi^{\Gamma}(X_t) ] =  \varphi^{\Gamma}(x) 
\quad \text{for all $ x \in S $ and $ t \ge 0 $} . 
\label{eq: Ginv}
\end{align}
\item[{\bf (A2)}] 
It holds that 
\begin{align}
\text{$ P_x[\Gamma_{\ve(q)}] \to 0 $ as $ q \down 0 $ for all $ x \in S^{\Gamma} $} , 
\label{}
\end{align}
where we abuse $ P_x $ for the extended probability measure of $ P_x $ 
supporting a standard exponential variable $ \ve $ independent of $ \cF_{\infty } $ and 
we set $ \ve(q) = \ve/q $ for $ q>0 $. 
\end{enumerate}

Note that, by the dominated convergence theorem, 
the condition {\bf (A2)} follows from the following condition: 
\begin{enumerate}
\item[({\bf A2$ ' $})] 
It holds that 
\begin{align}
\text{$ P_x[\Gamma_t] \to 0 $ as $ t \to \infty $ for all $ x \in S^{\Gamma} $} . 
\label{}
\end{align}
\end{enumerate}

By the multiplicativity, 
the condition \eqref{eq: Ginv} is equivalent to the condition that 
\begin{align}
\text{$ (\Gamma_t \varphi^{\Gamma}(X_t))_{t \ge 0} $ is a right-continuous 
$ ((\cG_t)_{t \ge 0},P_x) $-martingale for all $ x \in S $} 
\label{}
\end{align}
(for right-continuity, see, e.g., \cite[Theorem 5.8]{Get}). 
Under {\bf (A1)}, for $ x \in S^{\Gamma} $, 
we may define a probability measure $ P^{\Gamma}_x $ on $ (\bD,\cF_{\infty }) $, 
which we call the {\em penalised probability} of $ P_x $ for $ \Gamma $, 
by the following (see Section \ref{sec: ext}): 
\begin{align}
P^{\Gamma}_x|_{\cF_t} 
= \frac{\Gamma_t \varphi^{\Gamma}(X_t)}{\varphi^{\Gamma}(x)} \cdot P_x|_{\cF_t} 
\quad \text{for all $ t \ge 0 $} . 
\label{eq: PGamma}
\end{align}
It is then immediate that the {\em penalised process} 
$ \{ X,\cF_{\infty },(P^{\Gamma}_x)_{x \in S} \} $ 
is a Markov process with respect to $ (\cF_t)_{t \ge 0} $. 

We write $ \tend{P}{} $ for convergence in probability. 
In addition to {\bf (A1)} and {\bf (A2)}, we also introduce the following assumptions: 
\begin{enumerate}
\item[{\bf (A3)}] 
There is a non-negative finite $ \cF_{\infty } $-measurable functional 
$ \Gamma_{\infty } $ such that 
\begin{align}
\text{$ P^{\Gamma}_x \rbra{ \Gamma_t \tend{}{t \to \infty } \Gamma_{\infty } > 0 } = 1 $ 
for all $ x \in S^{\Gamma} $.} 
\label{}
\end{align}
\end{enumerate}

Note that in many examples 
we have {\bf (A3)} and $ P_x(\liminf_{t \to \infty } \Gamma_t = 0) = 1 $, 
which implies that the two measures $ P^{\Gamma}_x $ and $ P_x $ 
are mutually singular on $ \cF_{\infty } $. 

The following is a routine argument. 

\begin{Prop} \label{prop: penal meas}
Let $ \Gamma $ be a multiplicative weight. 
Then the following hold. 
\begin{enumerate}
\item 
Under {\bf (A1)}, it holds that 
\begin{align}
\text{$ P^{\Gamma}_x(\tau^{\Gamma} = \infty)=1 $ 
for all $ x \in S^{\Gamma} $}. 
\label{}
\end{align}
\item 
Under {\bf (A1)}, {\bf (A2)} and {\bf (A3)}, 
it holds that 
\begin{align}
P^{\Gamma}_x \rbra{ \varphi^{\Gamma}(X_t) \tend{}{t \to \infty } \infty } = 1 
\quad \text{for all $ x \in S^{\Gamma} $}. 
\label{eq: lim is infty}
\end{align}
\end{enumerate}
\end{Prop}

\Proof{
(i) 
We apply the optional stopping theorem to 
the $ ((\cG_t)_{t \ge 0},P_x) $-martingale 
$ M_t := \Gamma_t \varphi^{\Gamma}(X_t) / \varphi^{\Gamma}(x) $ 
(by {\bf (A1)}) 
to see that 
\begin{align}
P^{\Gamma}_x(\tau^{\Gamma}>t) 
=& P_x \sbra{ M_t ; \tau^{\Gamma}>t } 
\label{} \\
=& P_x \sbra{ M_{t \wedge \tau^{\Gamma}} } - P_x \sbra{ M_{t \wedge \tau^{\Gamma}} 
; \tau^{\Gamma} \le t } 
\label{} \\
=& P_x[M_0] - P_x \sbra{ M_{\tau^{\Gamma}} ; \tau^{\Gamma} \le t } = 1 , 
\label{}
\end{align}
which implies that $ P^{\Gamma}_x(\tau^{\Gamma}=\infty )=1 $. 

(ii) 
Let $ 0 \le s \le t < \infty $ and $ A_s \in \cF_s $. 
We then have 
\begin{align}
\begin{split}
P^{\Gamma}_x \sbra{ \frac{1}{\Gamma_t \varphi^{\Gamma}(X_t)} ; A_s } 
=& \frac{1}{\varphi^{\Gamma}(x)} P_x (A_s , \ \tau^{\Gamma}>t) 
\\
\le& \frac{1}{\varphi^{\Gamma}(x)} P_x (A_s , \ \tau^{\Gamma}>s) 
= P^{\Gamma}_x \sbra{ \frac{1}{\Gamma_s \varphi^{\Gamma}(X_s)} ; A_s } . 
\end{split}
\label{eq: ineq Qx 1/varphi 1}
\end{align}
This shows that $ N_t := 1/\{ \Gamma_t \varphi^{\Gamma}(X_t) \} $ 
is a non-negative $ P^{\Gamma}_x $-supermartingale with respect to 
the completed filtration $ (\bar{\cF}^{P^{\Gamma}_x}_t)_{t \ge 0} $ of $ (\cF_t)_{t \ge 0} $, 
and consequently it converges $ P^{\Gamma}_x $-a.s. as $ t \to \infty $ 
to some random variable $ N_{\infty } $. 
By {\bf (A3)}, we see that 
\begin{align}
\frac{1}{\varphi^{\Gamma}(X_t)} 
= \Gamma_t N_t 
\tend{}{t \to \infty } 
\Gamma_{\infty } N_{\infty } 
\quad \text{$ P^{\Gamma}_x $-a.s.}, 
\label{eq: lim is GYinfty }
\end{align}
which implies $ 1/\varphi^{\Gamma}(X_{\ve(q)}) \tend{P^{\Gamma}_x}{q \down 0} 
\Gamma_{\infty } N_{\infty } $. 
Using Fatou's lemma, 
we obtain 
\begin{align}
P^{\Gamma}_x \sbra{ \Gamma_{\infty } N_{\infty } } 
\le \liminf_{q \down 0} P^{\Gamma}_x \sbra{ \frac{1}{\varphi^{\Gamma}(X_{\ve(q)})} } 
= \frac{1}{\varphi^{\Gamma}(x)} \lim_{q \down 0 } P_x \sbra{ \Gamma_{\ve(q)} } = 0 
\label{eq: prf of varphiX_t to infty}
\end{align}
by {\bf (A2)}. Hence we obtain \eqref{eq: lim is infty}. 
}

\section{Subsequent Markov property} \label{sec: subs M} 

Let $ \Gamma $ be a multiplicative weight 
satisfying {\bf (A1)}, {\bf (A2)} and {\bf (A3)}. 
For $ x \in S^{\Gamma} $, we may define a measure $ \sP^{\Gamma}_x $ on $ (\bD,\cF_{\infty }) $, 
which we call the {\em unweighted measure} of $ P^{\Gamma}_x $, 
by 
\begin{align}
\sP^{\Gamma}_x = \varphi^{\Gamma}(x) \Gamma_{\infty }^{-1} \cdot P^{\Gamma}_x 
\quad \text{on $ \cF_{\infty } $}. 
\label{}
\end{align}
Note that $ \sP^{\Gamma}_x $ is $ \sigma $-finite on $ \cF_{\infty } $, 
because $ \bD = \bigcup_{n \in \bN} \{ \Gamma_{\infty } > 1/n \} $, $ \sP^{\Gamma}_x $-a.e. and 
\begin{align}
\sP^{\Gamma}_x(\Gamma_{\infty } > 1/n) \le n \varphi^{\Gamma}(x) < \infty 
\quad \text{for all $ n \in \bN $} . 
\label{}
\end{align}
The family of the unweighted measures satisfies the following property. 

\begin{Thm} \label{subsM}
Let $ \Gamma $ be a multiplicative weight 
satisfying {\bf (A1)}-{\bf (A3)}. 
Then, for any $ x \in S^{\Gamma} $, 
any non-negative $ \cF_t $-measurable functional $ F_t $ 
and any non-negative $ \cF_{\infty } $-measurable functional $ G $, 
it holds that 
\begin{align}
\sP^{\Gamma}_x \sbra{ F_t (G \circ \theta_t) } 
= P_x \sbra{ F_t \sP^{\Gamma}_{X_t}[G] ; \tau^{\Gamma}>t } . 
\label{eq: subs M}
\end{align}
\end{Thm}

\Proof{
By definition of $ \sP^{\Gamma}_x $, we have 
\begin{align}
\sP^{\Gamma}_x \sbra{ (F_t \Gamma_t) ((G \Gamma_{\infty }) \circ \theta_t) } 
=& \sP^{\Gamma}_x \sbra{ F_t (G \circ \theta_t) \Gamma_{\infty } } 
\label{} \\
=& \varphi^{\Gamma}(x) P^{\Gamma}_x \sbra{ F_t (G \circ \theta_t) } . 
\label{eq: subsequent M1}
\end{align}
By the Markov property for $ X $ under $ P^{\Gamma}_x $, 
by the local equivalence between $ P^{\Gamma}_x $ and $ P_x $, 
and by the global equivalence between $ P^{\Gamma}_x $ and $ \sP^{\Gamma}_x $, 
we obtain 
\begin{align}
\text{\eqref{eq: subsequent M1}} 
=& \varphi^{\Gamma}(x) P^{\Gamma}_x \sbra{ F_t P^{\Gamma}_{X_t}[G] } 
\label{} \\
=& P_x \sbra{ F_t \varphi^{\Gamma}(X_t) \Gamma_t P^{\Gamma}_{X_t}[G] } 
\label{} \\
=& P_x \sbra{ F_t \Gamma_t \sP_{X_t}[G \Gamma_{\infty }] } , 
\label{}
\end{align}
where we used the fact obtained from Proposition \ref{prop: penal meas} that 
$ X_t \in S^{\Gamma} $, $ P_x $-a.s. on $ \{ \Gamma_t > 0 \} $. 
Thus we obtain 
\begin{align}
\sP^{\Gamma}_x \sbra{ F_t \Gamma_t (G \Gamma_{\infty }) \circ \theta_t } 
=& P_x \sbra{ F_t \Gamma_t \sP^{\Gamma}_{X_t}[G \Gamma_{\infty }] } . 
\label{}
\end{align}
Replacing $ F_t $ by $ F_t \Gamma_t^{-1} 1_{\{ \tau^{\Gamma} > t \}} $ 
and $ G $ by $ G \Gamma_{\infty }^{-1} 1_{\{ \Gamma_{\infty } > 0 \}} $, 
we obtain the desired identity, 
since $ \tau^{\Gamma} = \infty $ and $ \Gamma_{\infty } > 0 $, $ \sP^{\Gamma}_x $-a.e. 
The proof is now complete. 
}

Theorem \ref{subsM} asserts that, 
the process under $ \sP^{\Gamma}_x $ 
behaves until a fixed time $ t $ 
as the process under $ P_x $ killed upon leaving $ S^{\Gamma} $, 
and it starts afresh at time $ t $ 
to behave as the process under $ \sP^{\Gamma}_{X_t} $. 
In this sense, we may call this property \eqref{eq: subs M} 
the {\em subsequent Markov property}.

\section{Universality class} \label{sec: univ}

Let $ \cE $ be a particular multiplicative weight satisfying {\bf (A1)}-{\bf (A3)}. 
We would like to give a sufficient condition 
for existence of a positive function $ c(x) $ such that 
\begin{align}
S^{\Gamma} \subset S^{\cE} 
\quad \text{and} \quad 
\sP^{\Gamma}_x = c(x) 1_{\{ \Gamma_{\infty } > 0 \}} \cdot \sP^{\cE}_x 
\quad \text{for all $ x \in S^{\Gamma} $}. 
\label{eq: sPGamma sPcE}
\end{align}
We note that [$ \sP^{\Gamma}_x = c(x) 1_{\{ \Gamma_{\infty } > 0 \}} \cdot \sP^{\cE}_x $] 
yields [$ \Gamma $ belongs to the universality class of $ \sP^{\cE}_x $] 
in the sense we mentioned in Introduction.

\begin{Thm}[Universality theorem] \label{thm: iden}
Let $ \cE $ and $ \Gamma $ be two multiplicative weights 
satisfying {\bf (A1)}-{\bf (A3)}. 
Suppose there exists a positive function $ c(x) $ such that 
\begin{align}
P^{\cE}_x \rbra{ \Gamma_t \tend{}{t \to \infty } \Gamma_{\infty } } = 1 
, \quad 
\frac{\varphi^{\Gamma}(X_t)}{\varphi^{\cE}(X_t)} \tend{P^{\cE}_x}{t \to \infty } c(x) 
\quad \text{for all $ x \in S^{\Gamma} $} 
\label{}
\end{align}
and 
\begin{align}
P^{\Gamma}_x \rbra{ \cE_t \tend{}{t \to \infty } \cE_{\infty } > 0 } = 1 
, \quad 
\frac{\varphi^{\Gamma}(X_t)}{\varphi^{\cE}(X_t)} \tend{P^{\Gamma}_x}{t \to \infty } c(x) 
\quad \text{for all $ x \in S^{\Gamma} $}. 
\label{}
\end{align}
(Notice that these assumptions do not follow from {\bf (A3)}.) 
Then \eqref{eq: sPGamma sPcE} holds. 
\end{Thm}

\Proof{
Let $ x \in S^{\Gamma} $ be fixed. 
Since $ P_x = P^{\Gamma}_x $ on $ \cF_0 $, we have 
\begin{align}
P_x(\cE_0=1) = P^{\Gamma}_x(\cE_0=1) \ge P^{\Gamma}_x(\cE_{\infty } > 0) = 1 , 
\label{}
\end{align}
which shows $ x \in S^{\cE} $. 
By the assumptions, we have 
\begin{align}
R_t \tend{P^{\cE}_x}{t \to \infty } R_{\infty } 
\ \text{and} \ 
R_t \tend{P^{\Gamma}_x}{t \to \infty } R_{\infty } 
\quad \text{with} \ 
R_t = \frac{\Gamma_t \varphi^{\Gamma}(X_t)}{\varphi^{\cE}(X_t)} 
\ \text{and} \ 
R_{\infty } = c(x) \Gamma_{\infty } . 
\label{eq: Iden 1}
\end{align}
Let $ s>0 $ 
and let $ F_s $ be a non-negative $ \cF_s $-measurable functional. 
For $ t>s $, we have 
\begin{align}
P^{\cE}_x \sbra{ F_s \cdot \frac{R_t}{1+R_t+\cE_t} } 
=& \frac{1}{\varphi^{\cE}(x) } 
P_x \sbra{ F_s \cdot \frac{R_t}{1+R_t+\cE_t} \cdot \cE_t \varphi^{\cE}(X_t) } 
\label{} \\
=& \frac{\varphi^{\Gamma}(x) }{\varphi^{\cE}(x) } 
P^{\Gamma}_x \sbra{ F_s \cdot \frac{R_t}{1+R_t+\cE_t} \cdot 
\frac{\cE_t \varphi^{\cE}(X_t)}{\Gamma_t \varphi^{\Gamma}(X_t)} } 
\label{} \\
=& \frac{\varphi^{\Gamma}(x) }{\varphi^{\cE}(x) } 
P^{\Gamma}_x \sbra{ F_s \cdot \frac{\cE_t}{1+R_t+\cE_t} } . 
\label{}
\end{align}
Letting $ t \to \infty $ and applying the dominated convergence theorem, we obtain 
\begin{align}
P^{\cE}_x \sbra{ F_s \cdot \frac{R_{\infty }}{1+R_{\infty }+\cE_{\infty }} } 
= \frac{\varphi^{\Gamma}(x) }{\varphi^{\cE}(x) } 
P^{\Gamma}_x \sbra{ F_s \cdot \frac{\cE_{\infty }}{1+R_{\infty }+\cE_{\infty }} } . 
\label{thm: iden2}
\end{align}
Since $ s>0 $ and $ F_s $ are arbitrary, we obtain 
\begin{align}
c(x) \Gamma_{\infty } \cdot P^{\cE}_x 
= \cE_{\infty } \cdot P^{\Gamma}_x , 
\label{}
\end{align}
which yields 
\begin{align}
c(x) 1_{\{ \Gamma_{\infty } > 0 \} } \cdot \sP^{\cE}_x 
= 1_{\{ \cE_{\infty } > 0 \} } \cdot \sP^{\Gamma}_x = \sP^{\Gamma}_x , 
\label{}
\end{align}
since $ P^{\Gamma}_x(\cE_{\infty }>0)=1 $. 
We thus obtain the desired result. 
}

\section{Penalisation problems} \label{sec: penal}

We give two systematic methods of ensuring the conditions {\bf (A1)} and {\bf (A2)} 
in penalisation problems. 

\subsection{Constant clock}

We give a general framework for 
penalisation problems with constant clock.

\begin{Prop} \label{penal lem1}
Let $ \Gamma $ be a multiplicative weight. 
Let $ \rho(t) $ be a function such that 
\begin{align}
\rho(t) \tend{}{t \to \infty } \infty 
\quad \text{and} \quad 
\frac{\rho(t)}{\rho(t-s)} \tend{}{t \to \infty } 1 
\quad \text{for all $ s>0 $}, 
\label{}
\end{align}
or in other words, $ r(\log t) $ is divergent and slowly varying at $ t = \infty $. 
Suppose there exists a process $ (M_s)_{s \ge 0} $ such that 
$ P_x(M_0 > 0) = P_x(\Gamma_0=1) $ for all $ x \in S $ 
and 
\begin{align}
\rho(t) P_x[ \Gamma_t | \cF_s ] 
\tend{}{t \to \infty } 
M_s 
\quad \text{in $ L^1(P_x) $ for all $ x \in S $ and all $ s \ge 0 $.} 
\label{eq: penal lem L1}
\end{align}
Then the weight $ \Gamma $ satisfies {\bf (A1)} and {\bf (A2$ ' $)} 
with 
\begin{align}
\varphi^{\Gamma}(x) = \lim_{t \to \infty } \rho(t) P_x[\Gamma_t] , 
\label{eq: penal lem L2}
\end{align}
and the following penalisation limit with constant clock holds: 
\begin{align}
\frac{\Gamma_t \cdot P_x}{P_x[\Gamma_t]} \tend{}{t \to \infty } P^{\Gamma}_x 
\quad \text{along $ (\cF_s)_{s \ge 0} $ for all $ x \in S^{\Gamma} $} . 
\label{eq: penal limit cc}
\end{align}
\end{Prop}

\Proof{
The convergence \eqref{eq: penal lem L1} for $ s=0 $ becomes 
\eqref{eq: penal lem L2}. 
By the multiplicativity $ \Gamma_t = \Gamma_s \cdot (\Gamma_{t-s} \circ \theta_s) $ 
and by the Markov property, we have 
\begin{align}
\rho(t) P_x[ \Gamma_t | \cF_s ] 
= \frac{\rho(t)}{\rho(t-s)} \Gamma_s \cdot \rho(t-s) P_{X_s}[\Gamma_{t-s}] 
\tend{}{t \to \infty } \Gamma_s \varphi^{\Gamma}(X_s) \quad \text{in $ P_x $-a.s.}. 
\label{}
\end{align}
which yields $ M_s = \Gamma_s \varphi^{\Gamma}(X_s) $. 
Hence we have 
\begin{align}
P_x[\Gamma_t \varphi^{\Gamma}(X_t)] 
= \lim_{u \to \infty } \rho(u) P_x[ P_x[\Gamma_u | \cF_t] ] 
= \lim_{u \to \infty } \rho(u) P_x[\Gamma_u] 
= \varphi^{\Gamma}(x) , 
\label{}
\end{align}
which shows that {\bf (A1)} is satisfied. 
As $ \rho(t) \to \infty $, we obtain {\bf (A2$ ' $)}. 
For $ s>0 $ and for a bounded $ \cF_s $-measurable functional $ F_s $, we obtain 
\begin{align}
\rho(t) P_x[F_s \Gamma_t] 
= P_x[F_s \rho(t) P_x[\Gamma_t|\cF_s] ] 
\tend{}{t \to \infty } 
P_x[F_s M_s] 
= \varphi^{\Gamma}(x) P^{\Gamma}_x[F_s] . 
\label{}
\end{align}
This shows \eqref{eq: penal limit cc}. 
}

\subsection{Exponential clock}

Conditioning and penalisation problems with exponential clock 
have been widely studied; 
see 
\cite{MR1419491}, \cite{MR2164035}, 
\cite{MR3689384}, \cite{MR3444297} and \cite{MR3909919}. 
We give a general framework for them.

\begin{Prop} \label{penal lem2} 
Let $ r(q) $ be a function defined for small $ q>0 $ such that 
$ r(q) \to \infty $ as $ q \down 0 $. 
We abuse $ P_x $ for the extended probability measure of $ P_x $ 
supporting a standard exponential variable $ \ve $ independent of $ (\cF_t)_{t \ge 0} $ and 
set $ \ve(q) = \ve/q $ for $ q>0 $. 
Suppose there exists a process $ (M_s)_{s \ge 0} $ such that 
$ P_x(M_0 > 0) = P_x(\Gamma_0=1) $ for all $ x \in S $ 
and 
\begin{align}
\begin{split}
\lim_{q \down 0} r(q) P_x[ \Gamma_{\ve(q)} | \cF_s ] 
= \lim_{q \down 0} r(q) P_x[ \Gamma_{\ve(q)} 1_{\{ \ve(q) > s \}} | \cF_s ] 
= M_s 
\quad \text{in $ L^1(P_x) $ } 
\\
\quad \text{for all $ x \in S $ and all $ s \ge 0 $.} 
\end{split}
\label{eq: penal lem2 L1}
\end{align}
Then the weight $ \Gamma $ satisfies {\bf (A1)} and {\bf (A2)} 
with 
\begin{align}
\varphi^{\Gamma}(x) = \lim_{q \down 0} r(q) P_x[\Gamma_{\ve(q)}] , 
\label{eq: penal lem2 L2}
\end{align}
and the following penalisation limit with exponential clock holds: 
\begin{align}
\lim_{q \down 0} \frac{\Gamma_{\ve(q)} \cdot P_x}{P_x[\Gamma_{\ve(q)}]} 
= \lim_{q \down 0} \frac{\Gamma_{\ve(q)} 1_{\{ \ve(q) > s \}} \cdot P_x}
{P_x[\Gamma_{\ve(q)} ; \ve(q) > s]} 
= P^{\Gamma}_x 
\quad \text{along $ (\cF_s)_{s \ge 0} $ for all $ x \in S^{\Gamma} $}. 
\label{eq: penal limit ec}
\end{align}
\end{Prop}

\Proof{
The convergence \eqref{eq: penal lem2 L1} for $ s=0 $ becomes \eqref{eq: penal lem2 L2}. 
By the multiplicativity $ \Gamma_t = \Gamma_s \cdot (\Gamma_{t-s} \circ \theta_s) $, 
by the Markov property and by the memoryless property 
\begin{align}
P(\ve(q)-s \in \cdot \mid \ve(q)>s) = P(\ve(q) \in \cdot) , 
\label{}
\end{align}
we have 
\begin{align}
r(q) P_x[ \Gamma_{\ve(q)} 1_{\{ \ve(q) > s \}} | \cF_s ] 
=& \e^{-qs} r(q) P_x[ \Gamma_{\ve(q)+s} | \cF_s ] 
\label{} \\
=& \e^{-qs} \Gamma_s r(q) P_{X_s}[\Gamma_{\ve(q)}] 
\tend{}{q \down 0} \Gamma_s \varphi^{\Gamma}(X_s) \quad \text{$ P_x $-a.s.}, 
\label{}
\end{align}
which yields $ M_s = \Gamma_s \varphi^{\Gamma}(X_s) $. 
Hence we obtain 
\begin{align}
P_x[\Gamma_t \varphi^{\Gamma}(X_t)] 
= P_x[M_t] 
= \lim_{q \down 0} r(q) P_x[ P_x[ \Gamma_{\ve(q)} | \cF_t ] ] 
= \lim_{q \down 0} r(q) P_x[ \Gamma_{\ve(q)} ] 
= \varphi^{\Gamma}(x) , 
\label{}
\end{align}
which shows that {\bf (A1)} is satisfied. 
As $ r(q) \to \infty $, we obtain {\bf (A2)}. 
For $ s>0 $ and for a bounded $ \cF_s $-measurable functional $ F_s $, we obtain 
\begin{align}
r(q) P_x[F_s \Gamma_{\ve(q)}] 
= P_x[F_s r(q) P_x[\Gamma_{\ve(q)}|\cF_s] ] 
\tend{}{q \down 0} 
P_x[F_s M_s] 
= \varphi^{\Gamma}(x) P^{\Gamma}_x[F_s] . 
\label{}
\end{align}
This shows \eqref{eq: penal limit ec}. 
}

\section{Brownian penalisation revisited} \label{sec: BM}

Let us look at some results of Roynette--Vallois--Yor \cite{MR2229621,MR2253307} 
and Najnudel--Roynette--Yor \cite{NRY} 
in our framework. 

Let $ \{ B=(B_t)_{t \ge 0},(W_x)_{x \in \bR} \} $ denote the canonical representation 
of the one-dimensional Brownian motion with $ W_x(B_0=x)=1 $. 
Set $ \bar{B}_t = \sup_{s \le t} B_s $ and let $ L_t $ denote the local time of $ B $ at 0. 
For the shift operator on the path space, we have 
\begin{align}
B_{t+s} = B_t \circ \theta_s 
, \quad 
\bar{B}_{t+s} = \bar{B}_s \vee (\bar{B}_t \circ \theta_s) 
, \quad 
L_{t+s} = L_s + (L_t \circ \theta_s) . 
\label{eq: shift B barB L}
\end{align}

For a technical reason, 
we set 
\begin{align}
S = \{ (x,y,l) \in \bR^3 : y \ge x , \ l \ge 0 \} 
\label{}
\end{align}
as the state space and consider the coordinate process 
$ X = (X_t)_{t \ge 0} = (X^B_t,X^{\rm sup}_t,X^{\rm lt}_t)_{t \ge 0} $ 
on the space of c\`adl\`ag paths from $ [0,\infty ) $ to $ S $. 
Writing $ a \vee b = \max \{ a,b \} $, 
we define $ P_{(x,y,l)} $ by the law on $ \bD $ of $ (B,y \vee \bar{B},l+L) $ under $ W_x $, 
and adopt the notation of Section \ref{sec: penal meas}. 
By the identities \eqref{eq: shift B barB L}, 
we see that the process $ \{ X,\cF_{\infty },(P_{(x,y,l)})_{(x,y,l) \in S} \} $ 
is a strong Markov process 
with respect to the augmented filtration. 

(1) Supremum penalisation. 
For an integrable function $ f:\bR \to [0,\infty ) $ 
such that for some $ - \infty < y_0 \le \infty $  
we have $ f(y)>0 $ for $ y \le y_0 $ and $ f(y) = 0 $ for $ y > y_0 $, 
we set 
\begin{align}
\Gamma^{{\rm sup},f}_t = \frac{f(X^{\rm sup}_t)}{f(X^{\rm sup}_0)} 
1_{\{ X^{\rm sup}_t \le y_0 \}} 
, \quad 
S^{{\rm sup},f} = \{ (x,y,l) \in S : y \le y_0 \} . 
\label{}
\end{align}
Then we see that $ \Gamma^{{\rm sup},f} $ is a multiplicative weight 
with $ S^{\Gamma^{{\rm sup},f}} = S^{{\rm sup},f} $ 
(in what follows we will omit similar remarks). 
By Roynette--Vallois--Yor \cite[Theorem 3.6]{MR2253307}, 
we see that all the assumptions of Proposition \ref{penal lem1} are satisfied 
with $ \rho(t) = \sqrt{\pi t/2} $ and 
\begin{align}
\varphi^{{\rm sup},f}(x,y,l) = y - x + \frac{1}{f(y)} \int_y^{y_0} f(u) \d u 
, \quad 
(x,y,l) \in S^{{\rm sup},f} , 
\label{}
\end{align}
so that {\bf (A1)} and {\bf (A2$ ' $)} are satisfied. 
By the discussion of Roynette--Vallois--Yor \cite[Subsection 1.4]{MR2253307}, 
we can derive that 
\begin{align}
P^{{\rm sup},f}_{(x,y,l)}(X^{\rm sup}_{\infty }>a) 
= \frac{ \int_a^{y_0} f(u) \d u }{ (y-x)f(y) + \int_y^{y_0} f(u) \d u } 
, \quad y \le a < \infty , 
\label{eq: RVY14}
\end{align}
and hence that [$ X^{\rm sup}_t = X^{\rm sup}_{\infty } $ for large $ t $] 
and [$ \Gamma^{{\rm sup},f}_t \to \Gamma^{{\rm sup},f}_{\infty } > 0 $] 
$ P^{{\rm sup},f}_{(x,y,l)} $-a.s., 
which shows {\bf (A3)}. 
By (ii) of Proposition \ref{prop: penal meas}, 
we obtain the following known results: 
\begin{align}
P^{{\rm sup},f}_{(x,y,l)} \rbra{ 
X^B_t \to - \infty 
, \ 
\frac{\varphi^{{\rm sup},f}(X_t)}{|X^B_t|} \to 1 } = 1 . 
\label{}
\end{align}

(2) Local time penalisation. 
For an integrable function $ f:[0,\infty ) \to [0,\infty ) $ 
such that for some $ 0 \le l_0 \le \infty $  
we have $ f(l)>0 $ for $ l \le l_0 $ and $ f(l) = 0 $ for $ l > l_0 $, 
we set 
\begin{align}
\Gamma^{{\rm lt},f}_t = \frac{f(X^{\rm lt}_t)}{f(X^{\rm lt}_0)} 
1_{\{ X^{\rm lt}_t \le l_0 \}} 
, \quad 
S^{{\rm lt},f} = \{ (x,y,l) \in S : l \le l_0 \} . 
\label{}
\end{align}
By Roynette--Vallois--Yor \cite[Theorem 3.13 and Lemma 3.15]{MR2253307}, 
we see that all the assumptions of Proposition \ref{penal lem1} are satisfied 
with $ \rho(t) = \sqrt{\pi t/2} $ and 
\begin{align}
\varphi^{{\rm lt},f}(x,y,l) = |x| + \frac{1}{f(l)} \int_l^{l_0} f(u) \d u 
, \quad 
(x,y,l) \in S^{{\rm lt},f} , 
\label{}
\end{align}
so that {\bf (A1)} and {\bf (A2$ ' $)} are satisfied. 
Moreover, {\bf (A3)} is also satisfied and 
\begin{align}
P^{{\rm lt},f}_{(x,y,l)} \rbra{ X^B_t \to \pm \infty } 
= \frac{x^{\pm} f(l) + \frac{1}{2} \int_l^{l_0} f(u) \d u }{|x|f(l) + \int_l^{l_0} f(u) \d u } 
\label{}
\end{align}
with $ x^{\pm} = \max \{ \pm x,0 \} $. 
It is then obvious that 
\begin{align}
P^{{\rm lt},f}_{(x,y,l)} \rbra{ 
|X^B_t| \to \infty 
, \ 
\frac{\varphi^{{\rm lt},f}(X_t)}{|X^B_t|} \to 1 } = 1 . 
\label{}
\end{align}

Note that the conditioning to avoid zero, 
which we have mentioned in Introduction, can be regarded as a special case 
of the local time penalisation with the weight 
$ 1_{\{ X^{\rm lt}_t = 0 \}} = \Gamma^{{\rm lt},f}_t $ for $ f(l) = 1_{\{ l=0 \}} $. 

(3) Kac killing penalisation with integrable potential. 
For an integrable function $ v:\bR \to [0,\infty ) $ satisfying 
\begin{align}
0 < \int_{\bR} (1+|x|) v(x) \d x < \infty , 
\label{}
\end{align}
we set 
\begin{align}
\Gamma^{{\rm Kac},v}_t = \exp \rbra{ - \int_0^t v(X^B_s) \d s } 
, \quad 
S^{{\rm Kac},v} = S . 
\label{}
\end{align}
By Roynette--Vallois--Yor \cite[Theorem 4.1]{MR2229621}, 
we see that all the assumptions of Proposition \ref{penal lem1} are satisfied 
with $ \rho(t) = \sqrt{\pi t/2} $ and 
$ \varphi^{{\rm Kac},v}(x,y,l) = \varphi_v(x) $ where $ \varphi_v $ 
is the unique solution to the Sturm--Liouville equation 
\begin{align}
\frac{1}{2} \frac{\d^2 \varphi_v}{\d x^2}(x) = v(x) \varphi_v(x) 
, \quad 
\lim_{x \to \pm \infty } \frac{\d \varphi_v}{\d x}(x) = \pm 1 . 
\label{eq: SLeq}
\end{align}
so that {\bf (A1)} and {\bf (A2$ ' $)} are satisfied. 
Moreover, {\bf (A3)} is also satisfied and 
\begin{align}
P^{{\rm Kac},v}_{(x,y,l)} \rbra{ X^B_t \to - \infty } 
= \frac{1}{C_v} \int_x^{\infty } \frac{\d y}{\varphi_v(y)^2} 
, \quad 
P^{{\rm Kac},v}_{(x,y,l)} \rbra{ X^B_t \to \infty } 
= \frac{1}{C_v} \int_{- \infty }^x \frac{\d y}{\varphi_v(y)^2} 
\label{}
\end{align}
with $ C_v = \int_{\bR} \frac{\d y}{\varphi_v(y)^2} $. 
By \eqref{eq: SLeq} it is obvious that 
\begin{align}
P^{{\rm Kac},v}_{(x,y,l)} \rbra{ 
|X^B_t| \to \infty 
, \ 
\frac{\varphi^{{\rm Kac},v}(X_t)}{|X^B_t|} \to 1 } = 1 . 
\label{}
\end{align}

(4) Kac killing penalisation with Heviside potential. 
For $ \lambda>0 $, we set 
\begin{align}
\Gamma^{{\rm Hev},\lambda}_t = \exp \rbra{ - \lambda \int_0^t 1_{ \{ X^B_s > 0 \} } \d s } 
, \quad 
S^{{\rm Kac},v} = S . 
\label{}
\end{align}
By Roynette--Vallois--Yor \cite[Theorem 5.1 and Example 5.4]{MR2229621}, 
we see that all the assumptions of Proposition \ref{penal lem1} are satisfied 
with $ \rho(t) = \sqrt{\pi t/2} $ and 
\begin{align}
\varphi^{{\rm Hev},\lambda}(x,y,l) = 
\begin{cases}
\frac{1}{\sqrt{2 \lambda}} \e^{- \sqrt{2 \lambda} x} & (x \ge 0) , \\
\frac{1}{\sqrt{2 \lambda}} - x & (x < 0) , 
\end{cases}
\label{}
\end{align}
so that {\bf (A1)} and {\bf (A2$ ' $)} are satisfied. 
Moreover, {\bf (A3)} is also satisfied and 
\begin{align}
P^{{\rm Hev},\lambda}_{(x,y,l)} \rbra{ 
X^B_t \to - \infty 
, \ 
\frac{\varphi^{{\rm Hev},\lambda}(X_t)}{|X^B_t|} \to 1 } = 1. 
\label{}
\end{align}

($*$) The universality class of Brownian penalisation. 
Take $ \cE_t = \exp (-X^{\rm lt}_t) $ as a special case of (2) with $ f(l) = \e^{-l} $. 
(Note that, by Najnudel--Roynette--Yor \cite[Theorem 1.1.2]{NRY}, 
the corresponding unweighted measure $ \sP^{\cE}_x $ 
coincides with $ \sW_x $ given in Introduction.) 
By the above argument, we see that 
all the assumptions of Theorem \ref{thm: iden} 
are satisfied with $ \cE $ and $ \Gamma = \Gamma^{{\rm sup},f} $, 
$ \Gamma^{{\rm lt},f} $, $ \Gamma^{{\rm Kac},v} $ or $ \Gamma^{{\rm Hev},\lambda} $, 
so that we obtain the following known result: 
\begin{align}
\sP^{\Gamma}_{(x,y,l)} = 1_{\{ \Gamma_{\infty }>0 \}} \cdot \sP^{\cE}_{(x,y,l)} 
\quad \text{for all $ (x,y,l) \in S^{\Gamma} $}. 
\label{}
\end{align}
We remark the following obvious facts: 
It holds up to $ \sP^{\cE}_{(x,y,l)} $-null sets that 
\begin{align}
\bD = \{ X^B_t \to \infty \ \text{or} \ X^B_t \to -\infty \} , 
\label{}
\end{align}
and that the event $ \{ \Gamma_{\infty }>0 \} $ becomes 
\begin{align}
\{ \Gamma^{{\rm sup},f}_{\infty }>0 \} 
=& \{ X^B_t \to - \infty \ \text{and} \ X^{\rm sup}_{\infty } \le y_0 \} , 
\label{} \\
\{ \Gamma^{{\rm lt},f}_{\infty }>0 \} 
=& \{ [X^B_t \to \infty \ \text{or} \ X^B_t \to -\infty] \ \text{and} \ 
X^{\rm lt}_{\infty } \le l_0 \} , 
\label{} \\
\{ \Gamma^{{\rm Kac},v}_{\infty }>0 \} 
=& \{ X^B_t \to \infty \ \text{or} \ X^B_t \to -\infty \} , 
\label{} \\
\{ \Gamma^{{\rm Hev},\lambda}_{\infty }>0 \} 
=& \{ X^B_t \to - \infty \} . 
\label{}
\end{align}

\section{L\'evy penalisation revisited} \label{sec: Levy} 

Let us look at some results of 
K.Yano--Y.Yano--Yor \cite{MR2552915,MR2744885}, Y.Yano \cite{MR3127911} 
and Takeda--K.Yano \cite{TY} 
in our framework. 

Let $ \{ Z=(Z_t)_{t \ge 0},(P^Z_x)_{x \in \bR} \} $ denote the canonical representation of 
one-dimensional strictly $ \alpha $-stable process of index $ 1 < \alpha < 2 $, 
skewness $ -1 \le \beta \le 1 $ and scaling parameter $ c_{\theta}>0 $: 
\begin{align}
P^Z_0[\e^{i \lambda Z_t}] =  \exp \rbra{ - c_{\theta} |\lambda|^{\alpha } 
\rbra{ 1 - i \beta \sgn (\lambda) \tan \frac{\pi \alpha }{2} } } 
, \quad \lambda \in \bR . 
\label{}
\end{align}
(For the facts in this paragraph, see e.g. \cite[Section VIII]{Ber}.) 
We assume that $ 1 < \alpha < 2 $ so as to exclude the Brownian case 
and to assure that zero is regular for itself: 
Writing $ T_0 = \inf \{ t>0 : Z_t=0 \} $ for the hitting time of zero, we have 
\begin{align}
P^Z_0(T_0 > 0) = 1 . 
\label{}
\end{align}
Set $ \bar{Z}_t = \sup_{s \le t} Z_s $ and let $ L_t $ denote the local time of $ Z $ at 0. 
Let 
\begin{align}
\rho := P^Z_0(Z_1>0) 
= \frac{1}{2} + \frac{1}{\pi \alpha } \arctan \rbra{ \beta \tan \frac{\pi \alpha }{2} } 
\in [1-1/\alpha ,1/\alpha ] 
\label{}
\end{align}
and let $ k $ denote the positive constant such that 
\begin{align}
\lim_{y \to \infty } y^{\alpha } P^Z_0(\bar{Z}>y) = k . 
\label{}
\end{align}

We set 
\begin{align}
S = \{ (x,y,l) \in \bR^3 : y \ge x , \ l \ge 0 \} 
\label{}
\end{align}
as the state space and consider the coordinate process 
$ X = (X_t)_{t \ge 0} = (X^Z_t,X^{\rm sup}_t,X^{\rm lt}_t)_{t \ge 0} $ 
on the space of c\`adl\`ag paths from $ [0,\infty ) $ to $ S $. 
We define $ P_{(x,y,l)} $ by the law on $ \bD $ of $ (Z,y \vee \bar{Z},l+L) $ under $ P^Z_x $, 
and adopt the notation of Section \ref{sec: penal meas}. 

(1) Supremum penalisation. 
For a non-increasing function $ f:\bR \to [0,\infty ) $ 
such that for some $ - \infty < y_0 \le \infty $  
we have $ f(y)>0 $ for $ y \le y_0 $ and $ f(y) = 0 $ for $ y > y_0 $, and 
\begin{align}
\int_0^{y_0} x^{\alpha \rho -1} f(y) \d y < \infty , 
\label{}
\end{align}
we set 
\begin{align}
\Gamma^{{\rm sup},f}_t = \frac{f(X^{\rm sup}_t)}{f(X^{\rm sup}_0)} 
1_{\{ X^{\rm sup}_t \le y_0 \}} 
, \quad 
S^{{\rm sup},f} = \{ (x,y,l) \in S : y \le y_0 \} . 
\label{}
\end{align}

By K.Yano--Y.Yano--Yor \cite[Theorem 5.1]{MR2744885}, 
we see that all the assumptions of Proposition \ref{penal lem1} are satisfied 
with $ \rho(t) = t^{\rho}/k $ and 
\begin{align}
\varphi^{{\rm sup},f}(x,y,l) = (y - x)^{\alpha \rho} 
+ \frac{\alpha \rho}{f(y)} \int_y^{y_0} f(u) (u-x)^{\alpha \rho-1} \d u 
, \quad 
(x,y,l) \in S^{{\rm sup},f} , 
\label{}
\end{align}
so that {\bf (A1)} and {\bf (A2$ ' $)} are satisfied. 
In the same way as that of deducing \eqref{eq: RVY14}, 
we see that [$ X^{\rm sup}_t = X^{\rm sup}_{\infty } $ for large $ t $] 
and [$ \Gamma^{{\rm sup},f}_t \to \Gamma^{{\rm sup},f}_{\infty } > 0 $] 
$ P^{{\rm sup},f}_{(x,y,l)} $-a.s., 
which shows {\bf (A3)}. 
By (ii) of Proposition \ref{prop: penal meas} 
and by the dominated convergence theorem, we obtain the following known results: 
\begin{align}
P^{{\rm sup},f}_{(x,y,l)} \rbra{ 
X^Z_t \to - \infty 
, \ 
\frac{\varphi^{{\rm sup},f}(X_t)}{(-X^Z_t)^{\alpha \rho}} \to 1 } = 1 . 
\label{eq: long time sup penal}
\end{align}

Note that the special case 
of the supremum penalisation with the weight 
$ 1_{\{ X^{\rm sup}_t = 0 \}} = \Gamma^{{\rm sup},f}_t $ for $ f(l) = 1_{\{ y=0 \}} $ 
corresponds to the conditioning to stay negative. 

(2) Local time penalisation. 
For an integrable function $ f:[0,\infty ) \to [0,\infty ) $ 
such that for some $ 0 \le l_0 \le \infty $  
we have $ f(l)>0 $ for $ l \le l_0 $ and $ f(l) = 0 $ for $ l > l_0 $, 
we set 
\begin{align}
\Gamma^{{\rm lt},f}_t = \frac{f(X^{\rm lt}_t)}{f(X^{\rm lt}_0)} 
1_{\{ X^{\rm lt}_t \le l_0 \}} 
, \quad 
S^{{\rm lt},f} = \{ (x,y,l) \in S : l \le l_0 \} . 
\label{}
\end{align}
By Takeda--K.Yano \cite{TY} and by certain computation in \cite[Section 5]{MR3072331}, 
we see that all the assumptions of Proposition \ref{penal lem2} are satisfied 
with $ r(t) = c_r q^{1/\alpha -1} $ for a certain constant $ c_r>0 $ and 
\begin{align}
\varphi^{{\rm lt},f}(x,y,l) = C_{\alpha ,\beta} (1-\beta \sgn(x)) |x|^{\alpha -1} 
+ \frac{1}{f(l)} \int_l^{l_0} f(u) \d u 
, \quad 
(x,y,l) \in S^{{\rm lt},f} 
\label{}
\end{align}
with a certain constant $ C_{\alpha ,\beta} > 0 $, 
so that {\bf (A1)} and {\bf (A2)} are satisfied. 
In the same way as that of deducing \eqref{eq: RVY14}, 
we see that [$ X^{\rm lt}_t = X^{\rm lt}_{\infty } $ for large $ t $] 
and [$ \Gamma^{{\rm lt},f}_t \to \Gamma^{{\rm lt},f}_{\infty } > 0 $] 
$ P^{{\rm lt},f}_{(x,y,l)} $-a.s., 
which shows {\bf (A3)}. 
By (ii) of Proposition \ref{prop: penal meas}, we obtain 
\begin{align}
P^{{\rm lt},f}_{(x,y,l)} \rbra{ 
(1-\beta \sgn(X^Z_t)) |X^Z_t|^{\alpha -1} \to \infty 
, \ 
\frac{\varphi^{{\rm lt},f}(X_t)}{C_{\alpha ,\beta} (1-\beta \sgn(X^Z_t)) |X^Z_t|^{\alpha -1}} 
\to 1 } = 1 ; 
\label{eq: long time lt penal2}
\end{align}
in particular, 
\begin{align}
P^{{\rm lt},f}_{(x,y,l)} \rbra{ 
X^Z_t \to - \infty 
, \ 
\frac{\varphi^{{\rm lt},f}(X_t)}{(-X^Z_t)^{\alpha -1}} \to 2 C_{\alpha ,1} } =& 1 
\quad (\text{if $ \beta = 1 $}), 
\label{eq: long time lt penal3} \\
P^{{\rm lt},f}_{(x,y,l)} \rbra{ 
X^Z_t \to \infty 
, \ 
\frac{\varphi^{{\rm lt},f}(X_t)}{(X^Z_t)^{\alpha -1}} \to 2 C_{\alpha ,-1} } =& 1 
\quad (\text{if $ \beta = -1 $}) . 
\label{eq: long time lt penal4}
\end{align}
In the case of $ -1 < \beta < 1 $, 
we have a stronger convergence result in Takeda--K.Yano \cite{TY}: 
\begin{align}
& P^{{\rm lt},f}_{(x,y,l)} \rbra{ \lim X^Z_t = \limsup X^Z_t = \limsup (-X^Z_t) = \infty } = 1 
\quad \text{if $ -1 < \beta < 1 $} . 
\label{eq: long time lt penal}
\end{align}

Note that the special case 
of the local time penalisation with the weight 
$ 1_{\{ X^{\rm lt}_t = 0 \}} = \Gamma^{{\rm lt},f}_t $ for $ f(l) = 1_{\{ l=0 \}} $ 
corresponds to the conditioning to avoid zero. 
See \cite{MR2648275} for comparison of two types of conditionings for L\'evy processes.

($*$) The universality classes of L\'evy penalisation. 
By \eqref{eq: long time sup penal}, 
it holds that 
\begin{align}
\{ \Gamma^{{\rm sup},f}_{\infty }>0 \} 
=& \{ X^Z_t \to -\infty \ \text{and} \ X^{\rm sup}_{\infty } \le y_0 \} 
\quad \text{up to $ \sP^{{\rm sup},f}_{(x,y,l)} $-null sets} 
\label{}
\end{align}
in any case of $ -1 \le \beta \le 1 $. 

($*$1) Consider the case of $ - 1 < \beta < 1 $. 
By \eqref{eq: long time lt penal}, 
it holds that 
\begin{align}
\begin{split}
\{ \Gamma^{{\rm lt},g}_{\infty }>0 \} 
=& \{ \lim X^Z_t = \limsup X^Z_t = \limsup (-X^Z_t) = \infty 
\ \text{and} \ X^{\rm lt}_{\infty } \le y_0 \} 
\\
& \hspace{16em} \text{up to $ \sP^{{\rm lt},g}_{(x,y,l)} $-null sets}. 
\end{split}
\label{}
\end{align}
This shows that the two $ \sigma $-finite measures 
$ \sP^{{\rm sup},f}_{(x,y,l)} $ and $ \sP^{{\rm lt},g}_{(x,y,l)} $ 
are singular to each other. 
Note that 
\eqref{eq: long time sup penal} and \eqref{eq: long time lt penal} imply 
\begin{align}
P^{{\rm sup},f}_{(x,y,l)} 
\rbra{ \frac{\varphi^{{\rm lt},g}(X_t)}{\varphi^{{\rm sup},f}(X_t)} \to 0 } = 1 
\label{}
\end{align}
because $ \alpha \rho > \alpha -1 $, 
so that the assumption of Theorem \ref{thm: iden} is not satisfied.

($*$2) Consider the case of $ \beta = 1 $, the spectrally positive case. 
Take $ \cE_t = \exp({X^{\rm sup}_0-X^{\rm sup}_t}) $ 
as a special case of (1) with $ f(y) = \e^{-y} $. 
Then, since $ \alpha \rho = \alpha -1 $, 
all the assumptions of Theorem \ref{thm: iden} 
are satisfied with $ \cE $ and 
$ \Gamma = \Gamma^{{\rm sup},f} $ or $ \Gamma^{{\rm lt},g} $, 
so that we conclude as a new result that 
\begin{align}
\sP^{\Gamma}_{(x,y,l)} = 1_{\{ \Gamma_{\infty }>0 \}} \cdot \sP^{\cE}_{(x,y,l)} 
\quad \text{for all $ (x,y,l) \in S^{\Gamma} $}. 
\label{}
\end{align}
It holds up to $ \sP^{\cE}_{(x,y,l)} $-null sets that 
\begin{align}
\bD = \{ X^Z_t \to - \infty \} , 
\label{}
\end{align}
and that the event $ \{ \Gamma_{\infty }>0 \} $ becomes 
\begin{align}
\{ \Gamma^{{\rm lt},g}_{\infty }>0 \} 
=& \{ X^Z_t \to - \infty \ \text{and} \ X^{\rm lt}_{\infty } \le l_0 \} . 
\label{}
\end{align}

($*$3) Consider the case of $ \beta = -1 $, the spectrally negative case. 
Then 
\begin{align}
\{ \Gamma^{{\rm lt},g}_{\infty }>0 \} 
= \{ X^Z_t \to \infty \ \text{and} \ X^{\rm lt}_{\infty } \le l_0 \} 
\quad \text{up to $ \sP^{{\rm lt},g}_{(x,y,l)} $-null sets}, 
\label{}
\end{align}
which shows that $ \sP^{{\rm sup},f}_{(x,y,l)} $ and $ \sP^{{\rm lt},g}_{(x,y,l)} $ 
are singular to each other.

\section{Langevin penalisation revisited} \label{sec: Lan}

Let us look at some results of Profeta \cite{MR3386368} 
in our framework. 

Let $ \{ (B,A),(W_{(b,a)})_{(b,a) \in \bR^2} \} $ 
denote the canonical representation 
of the two-dimensional diffusion $ (B,A)=(B_t,A_t)_{t \ge 0} $ 
where $ B $ is a Brownian motion starting from $ b $ and 
\begin{align}
A_t = a + \int_0^t B_u \d u . 
\label{}
\end{align}
This two-dimensional diffusion is a special case of the \emph{Langevin process} 
and the process $ A $ is called the \emph{integrated Brownian motion}. 
Set $ \bar{A}_t := \sup_{s \le t} A_s $. 

We set 
\begin{align}
S = \{ (b,a,y) \in \bR^3 : y \ge a \} 
\label{}
\end{align}
as the state space and consider the coordinate process 
\begin{align}
X = (X_t)_{t \ge 0} = (X^B_t,X^A_t,X^{\rm sup}_t)_{t \ge 0} 
\label{}
\end{align}
on the space of c\`adl\`ag paths from $ [0,\infty ) $ to $ S $. 
We define $ P_{(b,a,y)} $ by the law on $ \bD $ of 
$ (B,A,y \vee \bar{A}) $ under $ W_{(b,a)} $, 
and adopt the notation of Section \ref{sec: penal meas}. 

We recall the confluent hypergeometric function (see \cite[Chapter 13]{AI}): 
\begin{align}
U(\alpha, \beta, z) = \frac{1}{\Gamma(\alpha )} 
\int_0^{\infty } \e^{-zu} u^{\alpha -1} (1+u)^{\beta - \alpha -1} \d u 
, \quad \alpha >0 , \ \beta \in \bR , \ z>0 . 
\label{}
\end{align}
It is easy to see that 
\begin{align}
\frac{\d}{\d z} \rbra{ z^{\alpha } U(\alpha ,\beta,z) } 
= - \alpha (\beta - \alpha - 1) z^{\alpha -1} U(\alpha +1, \beta, z) . 
\label{eq: der U}
\end{align}
The following asymptotics are taken from \cite[Formulae 13.5.2 and 13.5.8]{AI}: 
\begin{align}
\lim_{z \to \infty } z^{\alpha } U(\alpha ,\beta,z) = 1 
\ (\beta \in \bR) 
, \quad 
\lim_{z \down 0} z^{\beta -1} U(\alpha ,\beta,z) = \frac{\Gamma(\beta-1)}{\Gamma(\alpha )} 
\ (1 < \beta < 2) . 
\label{eq: af}
\end{align}

(1) Conditioning to stay negative. 
We write $ \tau^A = \inf \{ t>0 : X^A_t \ge 0 \} $ 
for the exit time from $ (-\infty ,0) $ for the process $ X^A $ 
and set 
\begin{align}
\Gamma^A_t = 1_{\{ \tau^A > t \}} 
, \quad 
S^A = \{ (b,a,y) \in S : y < 0 \} 
= \{ (b,a,y) \in \bR^3 : a \le y < 0 \} . 
\label{}
\end{align}
By modifying Profeta \cite[Theorem 5]{MR3386368}, 
we see that all the assumptions of Proposition \ref{penal lem1} are satisfied 
with $ \rho(t) = c_1 t^{1/4} $ for a certain constant $ c_1>0 $ and 
\begin{align}
\varphi^A(b,a,y) = h(-a,-b) 
, \quad 
(b,a,y) \in S^A , 
\label{}
\end{align}
with a continuous function $ h:(0,\infty ) \times \bR \to (0,\infty ) $ given as 
\begin{align}
h(x,y) = 
\begin{cases}
(\frac92 x)^{1/6} z^{1/3} U(\frac16,\frac43,z) 
= y^{1/2} z^{1/6} U(\frac16,\frac43,z) 
& (y>0) , 
\\
\frac16 (\frac92 x)^{1/6} z^{1/3} U(\frac76,\frac43,z) \e^{-z} 
= \frac16 |y|^{1/2} z^{1/6} U(\frac76,\frac43,z) \e^{-z} 
& (y<0) , 
\end{cases}
\label{}
\end{align}
for $ x>0 $ and $ z = \frac{2}{9} \frac{|y|^3}{x} $, 
so that {\bf (A1)} and {\bf (A2$ ' $)} are satisfied. 
Moreover, {\bf (A3)} is also satisfied and 
\begin{align}
P^A_{(b,a,y)} \rbra{ 
X^B_t \to - \infty 
\ \text{and} \ 
X^A_t \to - \infty } = 1. 
\label{eq: XBXAinfty}
\end{align}
Let us prove this fact, 
as the part [$ X^B_t \to - \infty $] was not mentioned in \cite{MR3386368}. 
By the formulae \eqref{eq: af}, 
we see that both $ z^{1/6} U(\frac16,\frac43,z) $ and $ z^{1/6} U(\frac76,\frac43,z) \e^{-z} $ 
are bounded in $ z>0 $, we obtain $ h(x,y) \le c_2 |y|^{1/2} $ 
for some constant $ c_2>0 $. 
It holds $ P^A_{(b,a,y)} $-a.s. that, 
by (ii) of Proposition \ref{prop: penal meas}, 
\begin{align}
\varphi^A(X_t) = h(-X^A_t,-X^B_t) \to \infty , 
\label{cEXt h-XAt-XBt}
\end{align}
which yields [$ |X^B_t| \to \infty $]. 
But [$ P^A_{(b,a,y)}(X^B_t \to \infty) = 0 $], since 
[$ X^B_t \to \infty $] implies [$ X^A_t = a + \int_0^t X^B_s \d s \to \infty $], 
which contradicts the fact that $ X^A_0 = a < 0 $ and 
$ \tau^A = \infty $ by (i) of Proposition \ref{prop: penal meas}. 
Hence we obtain \eqref{eq: XBXAinfty}.

(2) Supremum penalisation. 
Let $ f:\bR \to [0,\infty ) $ be a continuous function 
such that 
for some $ - \infty < y_0 \le 0 $, 
we have $ f(y)>0 $ for $ y \le y_0 $ and $ f(y) = 0 $ for $ y > y_0 $. 
Set 
\begin{align}
\Gamma^{{\rm sup},f}_t = \frac{f(X^A_t)}{f(X^A_0)} 1_{\{ X^A_t \le y_0 \}} 
, \quad 
\text{
\begin{minipage}{16em}
$ S^{{\rm sup},f} = \{ (b,a,y) \in S : y \le y_0 \} $ 
\\
\hspace{2.35em} $ = \{ (b,a,y) \in \bR^3 : a \le y < y_0 \} $. 
\end{minipage}
}
\label{}
\end{align}
By Profeta \cite[Proposition 18 and Theorem 19]{MR3386368}, 
we see that all the assumptions of Proposition \ref{penal lem1} are satisfied 
with $ \rho(t) = c_1 t^{1/4} $ and 
\begin{align}
\varphi^{{\rm sup},f}(b,a,y) 
= h(y-a,-b) + \frac{1}{f(y)} \int_y^{y_0} f(w) \frac{\partial }{\partial w} h(w-a,-b) \d w 
, \quad 
(b,a,y) \in S^{{\rm sup},f} , 
\label{}
\end{align}
so that {\bf (A1)} and {\bf (A2$ ' $)} are satisfied. 
By a similar argument to that deducing \eqref{eq: RVY14}, 
we see that [$ X^{\rm sup}_t = X^{\rm sup}_{\infty } $ for large $ t $] 
$ P^{{\rm sup},f}_{(b,a,y)} $-a.s., 
and that [$ \Gamma^{{\rm sup},f}_t \to \Gamma^{{\rm sup},f}_{\infty } > 0 $] 
$ P^{{\rm sup},f}_{(b,a,y)} $-a.s., 
which shows {\bf (A3)}. 
By the fact that $ \frac{\partial h}{\partial w} \ge 0 $, we have 
\begin{align}
\varphi^{{\rm sup},f}(b,a,y) \le \rbra{ \sup_{y \le w \le y_0} f(w) } h(y_0-a,-b) . 
\label{}
\end{align}
By a similar argument after \eqref{cEXt h-XAt-XBt}, 
and by (ii) of Proposition \ref{prop: penal meas}, 
we can deduce 
\begin{align}
P^{{\rm sup},f}_{(b,a,y)} \rbra{ 
X^B_t \to - \infty 
\ \text{and} \ 
X^A_t \to - \infty } = 1 . 
\label{eq: XBXAinfty2}
\end{align}

($*$) The universality class of Langevin penalisation. 
We would like to compare the three unweighted measures 
$ \sP^A_{(b,a,y)} $, $ \sP^{{\rm sup},f}_{(b,a,y)} $ and $ \sP^B_{(b,a,y)} $. 
Here we write $ \tau^B = \inf \{ t>0 : X^B_t \ge 0 \} $ 
for the exit time from $ (-\infty ,0) $ for the Brownian motion $ X^B $ and set 
\begin{align}
\Gamma^B_t = 1_{\{ \tau^B > t \}} 
, \quad 
S^B = \{ (b,a,y) \in S : b < 0 \} . 
\label{}
\end{align}
The penalisation for the weight $ \Gamma^B $ 
is nothing else but the conditioning to stay negative for the Brownian motion, 
so that we obtain $ \varphi^B(b,a,y) = -b $. 
The penalized probability $ P^B_{(b,a,y)} $ is the minus times 3-dimensional Bessel process 
and the corresponding unweighted measure is given as 
$ \sP^B_{(b,a,y)} = (-b) P^B_{(b,a,y)} $. 
Since $ X^A_t = a + \int_0^t X^B_u \d u $, we obtain 
\begin{align}
P^B_{(b,a,y)} \rbra{ 
X^B_t \to - \infty 
\ \text{and} \ 
X^A_t \to - \infty } = 1 . 
\label{}
\end{align}

We prove the following proposition with conjectured assumptions. 

\begin{Prop}
Set $ \displaystyle Z_t = \frac{(-X^B_t)^3}{(-X^A_t)} $. 
Then the following assertions hold: 
\begin{enumerate}

\item 
Suppose the following conjecture is true: 
\begin{align}
Z_t \tend{P^A_{(b,a,y)}}{t \to \infty } \infty 
\ \text{and} \ 
Z_t \tend{P^{{\rm sup},f}_{(b,a,y)}}{t \to \infty } \infty 
\ \text{for} \ 
(b,a,y) \in S^{{\rm sup},f} . 
\label{}
\end{align}
Then $ \sP^{{\rm sup},f}_{(b,a,y)} $ and $ \sP^A_{(b,a,y)} $ coincide 
for $ (b,a,y) \in S^{{\rm sup},f} (\subset S^A) $. 

\item 
Suppose the following conjecture is true: 
\begin{align}
Z_t \tend{P^A_{(b,a,y)}}{t \to \infty } \infty 
\ \text{and} \ 
Z_t \tend{P^B_{(b,a,y)}}{t \to \infty } \infty 
\ \text{for} \ 
(b,a,y) \in S^A \cap S^B . 
\label{}
\end{align}
Then $ \sP^A_{(b,a,y)} $ and $ \sP^B_{(b,a,y)} $ 
are singular to each other for $ (b,a,y) \in S^A \cap S^B $. 

\end{enumerate}
\end{Prop}

\Proof{
(i) 
Set $ \displaystyle Z^{\rm sup}_t = \frac{(-X^B_t)^3}{(X^{\rm sup}_t-X^A_t)} $. 
Then $ Z_t \tend{P}{t \to \infty } \infty $ 
both for $ P = P^A_{(b,a,y)} $ and for $ P = P^{{\rm sup},f}_{(b,a,y)} $. 
Since $ X^B_t < 0 $ for large $ t $, we have 
\begin{align}
\frac{h(X^{\rm sup}_t - X^A_t,-X^B_t)}{h(- X^A_t,-X^B_t)} 
= \frac{(Z^{\rm sup}_t)^{1/6} U(\frac16,\frac43,Z^{\rm sup}_t)}
{(Z_t)^{1/6} U(\frac16,\frac43,Z_t)} 
\tend{P}{t \to \infty } 1 
\label{}
\end{align}
by the assumption. Noting that \eqref{eq: der U} implies 
\begin{align}
\frac{\partial }{\partial x} h(x,y) 
= c_3 x^{-5/6} \cdot z^{7/6} U({\textstyle \frac76,\frac43,z}) 
\le c_4 x^{-5/6} 
, \quad x,y>0 , \ z = \frac{2}{9} \frac{|y|^3}{x} 
\label{}
\end{align}
for some constants $ c_3,c_4 > 0 $, we obtain 
\begin{align}
\frac{\varphi^{{\rm sup},f}(X_t)}{\varphi^A(X_t)} 
\tend{P}{t \to \infty } 1 
\label{}
\end{align}
both for $ P = P^A_{(b,a,y)} $ and for $ P = P^{{\rm sup},f}_{(b,a,y)} $. 
We may now apply Theorem \ref{thm: iden} 
for $ \cE = \Gamma^A $ and $ \Gamma = \Gamma^{{\rm sup},f} $, 
and thus we obtain the desired result. 

(ii) 
By the assumption, we have 
\begin{align}
R_t := \frac{\Gamma^A_t \varphi^A(X_t)}{\varphi^B(X_t)} 
= \frac{\Gamma^A_t \cdot (-X^B_t)^{1/2} \cdot (Z_t)^{1/6} U(\frac16,\frac43,Z_t)}{(-X^B_t)} 
\tend{P}{t \to \infty } 0 
\label{}
\end{align}
both for $ P = P^A_{(b,a,y)} $ and for $ P = P^B_{(b,a,y)} $. 
By the same argument of Theorem \ref{thm: iden} 
with $ \cE = \Gamma^B $ and $ \Gamma = \Gamma^A $, we obtain 
\begin{align}
P^B_{(b,a,y)} \sbra{ F_s \cdot \frac{R_t}{1 + R_t + \Gamma^B_t} } 
= \frac{\varphi^A(b,a,y)}{\varphi^B(b,a,y)} 
P^A_{(b,a,y)} \sbra{ F_s \cdot \frac{\Gamma^B_t}{1 + R_t + \Gamma^B_t} } . 
\label{}
\end{align}
Letting $ t \to \infty $, we obtain $ P^A_{(b,a,y)}(\Gamma^B_{\infty }>0) = 0 $. 
Since $ P^B_{(b,a,y)}(\Gamma^B_{\infty }>0) = 1 $, we obtain the desired result. 
}

\section{Appendix: Extension of transformed probability measures} \label{sec: ext}

We discuss in general extension of the transformed probability measures 
given by local absolute continuity like \eqref{eq: PGamma}. 
Recall that $ \bD $ is the space of c\`adl\`ag paths 
from $ [0,\infty ) $ to a locally compact separable metric space $ S $ 
and $ X $ is the coordinate process on $ \bD $. 

\begin{Thm}
Let $ P $ be a probability measure on $ (\bD,\sigma(X)) $ 
and let $ (M_t)_{t \ge 0} $ be a non-negative martingale 
such that $ P[M_t]=1 $ for all $ t \ge 0 $. 
Then there exists a unique probability measure $ Q $ on $ (\bD,\sigma(X)) $ 
such that 
\begin{align}
Q|_{\cF^X_t} = M_t \cdot P|_{\cF^X_t} 
, \quad t \ge 0 , 
\label{}
\end{align}
where $ \cF^X_t = \sigma(X_s:s \le t) $ is the natural filtration of $ X $. 
\end{Thm}

\Proof{
Since $ \bigcup_{t \ge 0} \cF^X_t $ is a $ \pi $-system generating $ \sigma(X) $, 
uniqueness of $ Q $ follows immediately from Dynkin's $ \pi $-$ \lambda $ theorem. 

Let us prove existence of $ Q $. 
For $ n \in \bN $, let $ \bD_n $ denote 
the space of c\`adl\`ag paths from $ [n-1,n) $ to $ S $, 
equipped with the $ \sigma $-field 
$ \cB_n $ generated by the coordinate process on $ \bD_n $. 
We thus see that $ \bD $ is the product space of $ \{ \bD_n \} $: 
\begin{align}
\bD = \prod_{n=1}^{\infty } \bD_n 
, \quad 
\sigma(X) = \sigma \rbra{ \prod_{k=1}^n B_k \times \prod_{k=n+1}^{\infty } \bD_k : 
B_1 \in \cB_1, \ldots, B_n \in \cB_n ; \ n \in \bN } . 
\label{}
\end{align}

Let $ \mu_n $ denote the law on $ \bD_1 \times \cdots \times \bD_n $, 
the space of c\`adl\`ag paths from $ [0,n) $ to $ S $, 
of $ (X_t)_{0 \le t < n} $ under $ M_n \cdot P|_{\cF^X_n} $. 
We then see that $ \{ \mu_n \} $ is a projective sequence: 
\begin{align}
\mu_{n+1}(\cdot \times \bD_{n+1}) = \mu_n 
, \quad n \in \bN . 
\label{}
\end{align}
We may apply Daniell's extension theorem (cf. \cite[Theorem 6.14]{KalMod}) 
to see that there exists a sequence of random variables $ \{ \xi_n \} $ 
defined on a probability space $ (\Omega',\cF',P') $ 
such that $ \xi_n $ for each $ n $ takes values in $ \bD_n $ 
and the joint distribution of $ (\xi_1,\ldots,\xi_n) $ under $ P' $ for each $ n $ 
coincides with $ \mu_n $. 

We now define $ Q $ by the law on $ \bD $ 
of $ (\xi_1,\xi_2,\ldots) $ under $ P' $. 
For any $ A \in \cF^X_n $ for each $ n \in \bN $, 
we can find $ B \subset \bD_1 \times \cdots \times \bD_n $ which belongs to 
$ \sigma \rbra{ \prod_{k=1}^n B_k : B_1 \in \cB_1, \ldots, B_n \in \cB_n } $ 
such that $ A = \{ (X_t)_{0 \le t < n} \in B \} $, so that we obtain 
\begin{align}
Q(A) = P'((\xi_1,\ldots,\xi_n) \in B) 
= \mu_n(B) = P \sbra{ M_n ; (X_t)_{0 \le t < n} \in B } 
= P \sbra{ M_n ; A } . 
\label{}
\end{align}
We thus conclude that $ Q $ is as desired. 
}

\end{document}